\documentclass[a4paper]{article}                                                                  
\usepackage[a4paper]{geometry}  
\usepackage{graphicx}    
\usepackage{authblk}

\usepackage[hyphens]{url}
\usepackage{natbib} 
\usepackage{caption} 
\usepackage{amsmath, amssymb}
\usepackage{mathtools}
\usepackage{pgf}
\usepackage[capitalise]{cleveref}
\DeclareCaptionStyle{ruled}{labelfont=normalfont,labelsep=colon,strut=off} 
\usepackage{eurosym}

\usepackage{algorithm}
\usepackage{algorithmic}
\usepackage{booktabs}       
\usepackage{todonotes}
\usepackage{siunitx}
\usepackage{comment}

\usepackage{newfloat}
\usepackage{listings}
\usepackage{xspace}
\usepackage{subcaption}
\floatstyle{ruled}
\newfloat{listing}{tb}{lst}{}
\floatname{listing}{Listing}

\newcommand{\LR}{\textnormal{LR}}
\DeclareMathOperator*{\argmin}{arg\,min}

\usepackage{tabularx}

\title{Tricks from the Trade for Large-Scale Markdown Pricing: Heuristic Cut Generation for  Lagrangian Decomposition}

\begin{document}
\author[1]{Robert Streeck\footnote{robert.streeck@zalando.de}} 
\author[1]{Torsten Gellert} 
\author[1]{Andreas Schmitt} 
\author[1]{Asya Dipkaya} 
\author[1]{Vladimir~Fux} 
\author[1]{Tim Januschowski} 
\author[2]{Timo Berthold} 

\affil[1]{Zalando SE, Germany}                                              
\affil[2]{TU Berlin, Germany}                                                                                 

\date{}
\maketitle
\begin{abstract}
In automated decision making processes in the online fashion industry, the 'predict-then-optimize' paradigm is frequently applied, particularly for markdown pricing strategies. This typically involves a mixed-integer optimization step, which is crucial for maximizing profit and merchandise volume. In practice, the size and complexity of the optimization problem is prohibitive for using off-the-shelf solvers for mixed integer programs and specifically tailored approaches are a necessity. Our paper introduces specific heuristics designed to work alongside decomposition methods, leading to almost-optimal solutions. These heuristics, which include both primal heuristic methods and a cutting plane generation technique within a Lagrangian decomposition framework, are the core focus of the present paper. We provide empirical evidence for their effectiveness, drawing on real-world applications at Zalando SE, one of Europe's leading online fashion retailers, highlighting the practical value of our work. The contributions of this paper are deeply ingrained into  Zalando's production environment to its  large-scale catalog ranging in the millions of products and improving weekly profits by millions of Euros.
\end{abstract}

\section{Introduction}\label{sec:intro}

Online retailers face challenges that are especially important in periods where economies stagnate and companies focus on efficiency: they need to guarantee a smooth experience throughout the whole customer journey, organize distribution of large assortments (in the millions of products) and efficiently manage inventory across multiple countries. A crucial aspect for the success of an online retailer is its capability to dynamically scale commercial processes to account for sales periods, traffic growth and external market conditions. This puts strict constraints on the underlying algorithms and systems, which need to deliver results in tight timelines (minutes or hours), even if suboptimal, to ensure business continuity. The large assortment necessitates a high degree of automation so that human oversight cannot be guaranteed for both scale and response time requirements. Key operational questions are hence fully automated. 
A prime example for such a system that is both operationally and commercially critical and is subject to scaling and response-time requirements is markdown pricing.
 where 
The goal in markdown pricing is to maximize profit by managing a given inventory for the entire season.

In this paper, we consider large-scale price optimization systems for online retailers and the challenges that arise. Our show case throughout the paper will be Zalando SE, a leading online fashion retailer in Europe. With more than 50 million customers, presence in over 25 European countries and close to 2~million articles in its assortment, efficient pricing is one of the cornerstones of the Zalando business. To be successful, pricing systems and algorithms need to be 
\begin{enumerate}
    \item scalable, as the assortment is constantly growing to account for varying tastes and demands for novel propositions 
    \item fast, as the price update cycle has to meet commercial needs and match competitors' pacing.
\end{enumerate}
In addition, markdown pricing is one of the most effective levers to boost various dimensions of commercial performance (discounts depth, revenue, profit, etc.), both in specific countries/markets and on a company-wide level. This requires the pricing system to be able to steer towards strategic company goals, which in turn necessitates an integrated optimization (due to linking constraints) for all, or at least larger sets of markets. This poses an additional challenge for algorithm design, since the aforementioned commercial needs demand timely execution. 

A typical approach in markdown pricing follows the predict-then-optimize paradigm (see e.g.,~\citep{predict-then-optimize,prometheus,kedia2020price,caro12} for examples in pricing) where a number of forecasting models (e.g.,~\citep{kunz2023deep,schultz2023causal} for forecasting for pricing at Zalando and~\citep{salinas2020deepar,lim2021,rasul2021autoregressive} for more general examples) provide input to an optimization problem. In practice, mixed-integer programming (MIP) solutions are often employed given their adaptability to changing constraints~\citep{kedia2020price,Li2022,caro12}. The pricing approach of our show case Zalando follows this paradigm and here, we will consider in particular the "optimize" part. Previous work~\citep{Li2022} described the basic setup of the pricing optimization problem at Zalando and key ingredients necessary for its solution, such as a Lagrangian decomposition framework, modelling approaches and basic heuristics. However, in the years since the deployment of this solution in production, roughly since 2021, there were occasional repercussions from unfavourable practical situations that led to non-convergence of the solution. As a further motivation for improving the heuristic performance of the markdown-pricing system, the scale of large online retailers like Zalando implies that even modest percentual gains in finding more profitable prices result in large absolute monetary gains -- easily in the multi-digit million Euro range. 

The main contribution of this paper are extensions of the pricing algorithm described by Li et al.\ in~\citep{Li2022} via novel combinations of existing heuristic techniques both for the Cutting Plane iterations and primal solution construction. The focus is on providing high quality, but potentially sub-optimal solutions within tight time boundaries as commonly dictated by commercial processes. These extensions address the aforementioned challenges around scalability and speed. We also show in empirical evaluations that our approach allows to improve revenue by 3\%--6\% and profit by 2\%--5\% for our show case Zalando~SE while having negligible impact on solution speed.

Our paper is structured as follows. In~\cref{sec:rel}, we discuss work related to this paper. In~\cref{sec:background}, we introduce the necessary notation and give the formalization of the markdown pricing problem. \cref{sec:overview} discusses the algorithms used to solve the markdown pricing problem and in~\cref{sec:heuristics}, we present the main methodological novelties of this paper, namely the heuristics for solving the markdown problem. \cref{sec:exp} provides extensive experimental results. We conclude in~\cref{sec:concl}.

\section{Related Work}\label{sec:rel}

Pricing and revenue management are mature areas of research stretching multiple disciplines, see e.g.,~\citep{phillips2021pricing} for an introduction to the topic from an Economics perspective. The Machine Learning community typically covers aspects involving prediction, like demand forecasting, e.g.,~\citep{eisenach2020mqtransformer,kunz2023deep,oreshkin2020nbeats} and sometimes also taking a reinforcement learning perspective (e.g.,~\citep{madeka2022deep} for an example of an industrial application). Econometricians tend to focus on causal aspects or price elasticity estimation~\citep{athey07,STRAUSS2018375,TWFE,deaton1980economics,fogarty2010demand,hughes2008evidence,defusco2017interest} and in the formulation of the optimization problems needed for profit-optimal pricing~\citep{phillips2021pricing}. Here, we take an Operations Research perspective (see~\citep{predict-then-optimize,Li2022} for examples), where we take both the predictive problems and the problem formulation as given and rather work on solving the optimization problem most efficiently within the requirements dictated by production purposes.  This is, to the best of our knowledge, under-explored in this specific form in pricing.~\citet{CHEN2024718} present a notable exception with a recent contribution, however, not for markdown pricing but for multi-product pricing in a ride-hauling setting. 

Beyond pricing, in the larger context of mathematical optimization, heuristic methods are omnipresent, see~\citep{Berthold2014} for an overview. For complex MIP problems, a plethora of general-purpose heuristics exist, see, e.g.,~\citep{Berthold2014b,berthold2019,FischettiLodi03,GamrathBertholdHeinzetal2019}. Also, cut generation, aggregation, and selection strategies are often of a heuristic nature, consider, e.g., \citep{andersen2005reduce,BonamiCornuejolsDashEtAl2008Projected,Caprara1996,turner2023cutting} and in particular,~\citep{fischetti2011relax} which heuristically generates cutting planes from a Lagrangian relaxation of a general MIP. 
While our approach is general in the sense that it can be easily adapted to other applications of using Lagrangian relaxations to solve MIPs, we believe that the particular tuning of thresholds, limits, and selection criteria would not generalize as easily.

\section{Background and Problem Formulation}\label{sec:background}

    The pricing problem we tackle can be formalized as follows.
    For an assortment of $n$ articles, we want to determine for each article~$i \in N =\{1, \ldots, n\}$ a price~$x_i$ to maximize the long-term profit (LTP) $\max_{x_i \in \mathcal{X}_i} f_i(x_i)$.
    Here, the set~$\mathcal{X}_i$, with liberty in the notation, denotes the feasible set of a mixed integer programming formulation for the article prices, described by (linear) constraints, bounds, and auxiliary variables not explicitly mentioned.
    Similarly, the linear objective function $f(x)$ depends on auxiliary variables. The notation is deliberately kept general to demonstrate the broad possible applications.
    For the showcase of Zalando, this entails that~$x_i$ models varying prices for the article in different countries and various (future) weeks.
    Auxiliary variables model, amongst others, the sales and the stock of the article over the future weeks. In this showcase, we consider around 500\;000 articles. For a single article, the model contains up to 10\;000 variables (5\;000 binary) and 10\;000 constraints.
    More details on the formulation of the pricing problem at Zalando can be found in~\citep{Li2022}.
    
    Additional business constraints restrict the solution space. Typical examples are targets for an average discount rate or a gross merchandising volume, both aggregated over each country.
    These targets are often referred to as \emph{linking constraints}, as they depend on all articles.
    In a general form, the primal pricing problem is thus given as
    \begin{equation}
    P = \max_{\substack{x \in \mathcal{X}\\ A x \le b}} \sum_{i \in N} f_i(x_i), \label{eq:pricing_problem}
    \end{equation}
    where~$Ax\leq b$ enforces the linking constraints, the set~$\mathcal{X} = \mathcal{X}_1 \times \dots \times \mathcal{X}_n$ describes the Cartesian product of the individual articles feasible sets and~$x = (x_1, \dots, x_n)$.

    We again refer to~\citep{Li2022} for a comprehensive descriptions of the linking constraints and the basic Lagrangian algorithm which we summarize in the next section.
    For a concrete example of linking constraints, consider the average discount rate over all articles per given country which we typically constrain to a certain range as a leverage for business leaders to steer prices. The \emph{sales weighted discount rate} extends the discount rate by weighing discounts by sales, \[\textrm{sDR} = \frac{\sum_i ({p_i}-\bar{x}_i)\cdot s_{i,\bar{x}_i}}{\sum_i {p_i} \cdot s_{i,\bar{x}_i}}\] where $p_i$ denotes the undiscounted price of article~$i$ and $s_{i,\bar{x}_i}$ the number of items sold of article~$i$ when offered for price~$\bar{x}_i$, and the set of variables $\bar{x}_i$ is a subset of the decision variables~$x_i$.

    \subsubsection*{Lagrangian Relaxation}
    With typical sizes of online retailer assortments ranging in the millions of articles (which would result in billions of variables in our problem), the pricing problem is prohibitively large for off-the-shelf solvers. Techniques such as decomposition or aggregation become inevitable for identifying high-quality solutions in practice.
    
    A standard technique for decomposition approaches is based on a Lagrangian relaxation. Here, the violation of each linking constraint is penalized using Lagrangian multipliers~$\lambda\ge 0$ in the objective 
    \[
    \LR(\lambda, x) = \sum_{i \in N} f_i(x_i)  - \lambda^\top (b - A x).
    \] 
    For this modified objective the Lagrangian relaxation searches for an optimal solution~$x$ for given~$\lambda$, leading to the optimization problem 
     \[
     \LR(\lambda) = \max_{x \in \mathcal{X}} \LR(\lambda, x) = \sum_{i \in N} \underbrace{\max_{x_i\in \mathcal{X}_i} (f_i(x_i) + \lambda^\top A_i x_i)}_{\text{single article problem}} - \underbrace{\lambda^\top b}_{\text{constant offset}}
    \]
    
    By relaxing the linking constraints, we end up with one independent optimization problem per article as presented at the beginning of this section.
    An interesting aspect in solving the optimization problem comes via business requirements. For example, throughout a week, Zalando needs to re-optimize prices several times for its entire assortment of several hundred thousands of articles. The window for taking the markdown pricing decision is determined by the time when the most recent data becomes available and commercial steering meetings. This is typically a few hours. 
    In practice, a single article's optimal long-term profit can be computed well below 10~seconds on a single machine. 
    Via parallelization, the entire assortment can be solved in acceptable time (not exceeding a few hours), when enough computational resources are provided, i.e., a cluster of machines which provides around 1000~cores, provisioned by a cloud provider.
    The memory required for optimization of the entire assortment easily extends several hundreds of gigabyte.
    
    \section{Overview of Algorithms}
    \label{sec:overview}
    The objective value of the Lagrangian relaxation yields an upper bound on the pricing problem for every $\lambda \ge 0$.
    The optimization challenge is to find a combination of multipliers that provides a smallest upper bound, via the Lagrangian multiplier problem  $\mu^*=\min_{\lambda \geq 0} \LR(\lambda)$.
    We will describe how we managed to find suitable multipliers for the markdown pricing problem in the following.
    
    \subsubsection*{Cutting Plane Procedure}
    There are different approaches to obtain~$\mu^*$, e.g., Subgradient or Bundle methods, see \cite{guignard2003} for an overview. Here, we follow the approach by~\citep{li2021large} and use a Cutting Plane procedure. To motivate this approach, note that the multiplier problem to compute~$\mu^\ast$ can be reformulated as a linear problem where each constraint corresponds to a feasible solution $x \in \mathcal{X}$ of the pricing problem as
    \begin{equation}
        \mu^* = \min_{\lambda \ge 0, \mu} \mu ~\text{s.t.}~ \mu - (b- A x)^\top \lambda\ge \sum_{i \in N} f_i(x_i)  ~ \forall x \in \mathcal{X}. \label{prob:cp}
    \end{equation}
     Since every point in~$\mathcal{X}$ is a mixed-integer solution to a linear system of inequalities, it is computationally hard to obtain all constraints and therefore the model explicitly: the number of solutions can be exponential in the number of variables.

     The Cutting Plane approach avoids this problem by considering a relaxation to Problem~\eqref{prob:cp}, which contains only a subset of these potentially exponentially many constraints. The relaxation is tightened in every iteration by adding further valid inequalities.
     In more detail, assume that in iteration~$j$ a number of inequalities determined by points $X^1, \dots, X^j \in \mathcal{X}$ are included in the relaxation. Then the relaxation is given by
    \begin{equation}
        \mu^j, \lambda^j  = \argmin_{\lambda \ge 0, \mu} \mu ~\text{s.t.}~ \mu \ge \sum_{i \in N} f_i(X^k_i) + \lambda^\top(b - A X^k) ~ \forall k \le j. \label{prob:cp_approx}
    \end{equation}

    The solution to~\eqref{prob:cp_approx} yields a lower bound on $\mu^*$ and is simple to solve, see below.
    To check whether the relaxation solution solves the original problem, i.e., whether $\mu^j$, $\lambda^j$ belong to the feasible set of Problem~\eqref{prob:cp} it suffices to search for a point $x \in \mathcal{X}$ such that the respective inequality
    \begin{equation}
       \mu - (b - A x)^\top \lambda\ge \sum_{i \in N} f_i(x_i)\label{ineq:cut}
    \end{equation}
    is violated for $\mu^j, \lambda^j$. The separation problem is equivalent to computing $\LR(\lambda^j)$ and comparing its value against $\mu^j$. If $\LR(\lambda^j) > \mu^j$, the maximizing solution of $\LR(\lambda^j)$ is added as $X^{j+1}$. In the next iteration this leads to another inequality in the relaxation cutting off the current solution $\mu^j, \lambda^j$. Otherwise, the optimal Lagrangian solution value is found. 
    
    This approach enables us to evaluate the quality of our relaxation in each iteration. On the one hand, the value $\min_{k \le j} \LR(\lambda^k)$ yields an upper bound on $\mu^*$ and also $P$ and is therefore called the \emph{dual bound}. On the other hand, the \emph{relaxed primal bound}~$\mu^j$ yields a lower bound to $\mu^*$, as a relaxation is solved. 
 
    It is worth mentioning that the relaxation of the Cutting Plane problem~\eqref{prob:cp_approx} itself is very easy to solve. It is a linear problem with one variable for each linking constraint and with $j$~constraints.
     In contrast, computing violated cuts involves higher computational effort, since a single new cut requires all subproblems of the entire assortment to be solved. Therefore, it's crucial to aim for a small number of iterations to reach convergence.

    \subsubsection*{Practical considerations for the Cutting Plane procedure}

    To avoid that the relaxations given by~\eqref{prob:cp_approx} are unbounded during the first iterations, an upper bound $\lambda \le \bar{\lambda}$ is added to the linear problem. This upper bound is also useful in case a linking constraint is infeasible. The corresponding multiplier would otherwise grow indefinitely. We choose the upper bound~$\bar{\lambda}$ large enough that nearly all feasible linking constraints should have a non-violating solution to the Lagrangian relaxation with $\bar{\lambda}$. Consequently, we start the algorithm with $\lambda^0 = 0$, and it will choose $\lambda^1_\ell = \bar{\lambda}$ for all linking constraints~$\ell$ violated by $X^0$. In cases with a large number of linking constraints (and consequently a high-dimensional Cutting Plane problem) it might require a large number of iterations and added cuts until all occurrences of $\lambda^j_\ell = \bar{\lambda}$ for at least one~$\ell$ have been dropped from the cutting plane problem. Indeed, in practice, we observe that most early iterations will contain multipliers at $\bar{\lambda}$, even though such multipliers over-fulfill on the constraints and are likely too high.

    \subsubsection*{Primal Heuristic}

    The Cutting Plane procedure yields in each iteration a potential solution~$X^k$ with price suggestions for the assortment. However, even if the procedure is converged and the optimal value of the Lagrangian Relaxation is found, in practice the solutions commonly violate the linking constraints. This problem is amplified by the fact that there is not sufficient time to perform all iterations necessary for convergence. 

    To obtain a good solution from the information collected during the Cutting Plane procedure we use a MIP in order to provide a solution close to the optimum with only a minimal amount of constraint violations.

    Given the solutions~$X^k$ for $1\leq k \leq j$, the MIP selects for each article $i$ an iteration $k_i \le j$ such that the offer $x_i$  given by $x_i = X^{k_i}_i$ minimizes a combination of linking constraint violation and profit loss scaled by the highest violation per linking constraints and highest profit observed during the Cutting Plane procedure. The resulting MIP is then
          \begin{equation}
      \begin{aligned}
        \max ~& \bar{p} f(x) + \bar{v}^\top \delta& \\
         \mbox{s.t.}~&\delta \le Ax - b  \\
                     & x_i =  \sum_{1 \le k \le j} X^k_i {y_{i k}} &\forall 1 \le i \le n\\
                     & \sum_{1 \le k \le j} {y_{i k}} = 1 & \forall 1 \le i \le n\\
                     & {y_{i k}}\in \{0,1\}, {\delta} \le 0,        
      \end{aligned}
      \label{prob:primal_heuristic}
      \end{equation}
     where the binary variables $y_{i k}$ are one if and only if article $i$ takes the values of iteration $k$ and each $\delta_\ell$ measures the violation of the linking constraint~$\ell$. The constants $\frac{1}{\bar{p}} = \max_{k} f(X^k)$ and $\frac{1}{\bar{v}_\ell} = \max_{k} A_{\ell}^\top X^k - b_\ell - \min_{k} A_{\ell}^\top X^k - b_\ell$ for each linking constraint~$\ell$ are furthermore used as scaling factors.
    This MIP is easier to solve than the original primal problem~\eqref{eq:pricing_problem}, since the interaction of variables for one article described by~$\mathcal{X}_i$ and $f_i$ is already evaluated.
    We refer to~\cite{Li2022} for an in-depth explanation of this final feasibility step.

\section{Heuristics for the Markdown Pricing Problem}\label{sec:heuristics}
    
 \begin{figure}[t] 
     \begin{subfigure}{0.48\textwidth}
        \centering
        \begin{tikzpicture}[scale=0.9]
          \draw[thick,->] (-0.125,0) -- (5,0) node [right]{$\lambda$};
	   \draw[thick,->] (0,-0.125) -- (0,4.1) node [above]{$\mu, \textrm{LR}(\lambda)$};

        %draw unkown area
        \path[fill=green!80!black] (0,4)--(1,3)--(2,2.5)--(4,2.75)--(4.5,3)--(5,4)--cycle;

        %cuts existing
        \draw[-, blue!70!black, dashed] (-0.1,4.1) -- (0,4) -- (4,0) node [left, pos=0.25] {$\textnormal{cut}_1$} -- (4.1,-0.1);
        \draw[-, blue!70!black,dashed] (2.9,-0.2)--(3,0) -- (5,4) node [right, pos=0.95]{$\textnormal{cut}_2$} --(5.1,4.2) ;

        %draw points
        \filldraw (0.25,3.75) circle (1.5pt) node [below] {$X^1$};
        \draw [-] (0.25,0.1) -- (0.25,-0.1) node [below] {$\lambda^1$};

        \filldraw (4.6, 3.2) circle (1.5pt) node [right] {$X^2$};
        \draw [-] (4.6,0.1) -- (4.6,-0.1) node [below] {$\lambda^2$};

        \filldraw (10/3,2/3) circle (1.5pt) node [left] {$\mu^2$};
        \draw [-] (10/3,0.1) -- (10/3,-0.1) node [below] {$\lambda^3$};

        % additional cut 1
        %\draw[-, blue, dashed] (-0.2,0.7) -- (5.2,2.8) node [below]{$\textnormal{cut}_3$};
        %\filldraw (2.33,1.68) circle (1.5pt) node [right] {$z''$};
        %\draw [-] (2.33,0.1) -- (2.33,-0.1) node [below] {$\lambda''$};

         % additional cut 2
        %\draw[-, blue, dashed] (-0.2,2) -- (5.2,1.85) node [below]{$\textnormal{cut}_4$};
        %\filldraw (2.9,1.9) circle (1.5pt) node [above] {$z'''$};
        %\draw [-] (2.9,0.1) -- (2.9,-0.1) node [below] {$\lambda'''$};
        
         \end{tikzpicture} 
        \caption{With two solution for multiplier $\lambda^1$ and~$\lambda^2$ the Cutting Plane procedure chooses multiplier~$\lambda^3$ for the next iteration.}
    \end{subfigure}\hspace{\fill}
    \begin{subfigure}{0.48\textwidth}
        \centering
        \begin{tikzpicture}[scale=0.9]
                  \draw[thick,->] (-0.125,0) -- (5,0) node [right]{$\lambda$};
	   \draw[thick,->] (0,-0.125) -- (0,4.1) node [above]{$\mu, \textrm{LR}(\lambda)$};

        %draw unkown area
        \path[fill=green!80!black] (0,4)--(1,3)--(2,2.5)--(4,2.75)--(4.5,3)--(5,4)--cycle;

        %cuts existing
        \draw[-, blue!70!black, dashed] (-0.1,4.1) -- (0,4) -- (4,0) node [left, pos=0.25] {$\textnormal{cut}_1$} -- (4.1,-0.1);
        \draw[-, blue!70!black,dashed] (2.9,-0.2)--(3,0) -- (5,4) node [right, pos=0.95]{$\textnormal{cut}_2$} --(5.1,4.2) ;

        %draw points
        \filldraw (0.25,3.75) circle (1.5pt) node [below] {$X^1$};
        \draw [-] (0.25,0.1) -- (0.25,-0.1) node [below] {$\lambda^1$};

        \filldraw (4.6, 3.2) circle (1.5pt) node [right] {$X^2$};
        \draw [-] (4.6,0.1) -- (4.6,-0.1) node [below] {$\lambda^2$};

        \filldraw (10/3,2/3) circle (1.5pt) node [left] {$\mu^2$};
        %\draw [-] (10/3,0.1) -- (10/3,-0.1) node [below] {$\lambda_3$};

        % additional cut 1
        \draw[-, blue, dashed] (-0.2,0.7) -- (5.2,2.8) node [pos=0.8]{$\textnormal{heuristic  cut}$};
        \filldraw (2.33,1.68) circle (1.5pt) node [right] {$\mu^3$};
        \draw [-] (2.33,0.1) -- (2.33,-0.1) node [below] {$\lambda^3$};

         % additional cut 2
        %\draw[-, blue, dashed] (-0.2,2) -- (5.2,1.85) node [below]{$\textnormal{cut}_4$};
        %\filldraw (2.9,1.9) circle (1.5pt) node [above] {$z'''$};
        %\draw [-] (2.9,0.1) -- (2.9,-0.1) node [below] {$\lambda'''$};
         \end{tikzpicture}
        \caption{An additional cut improves the bound $\mu^2$ provided by the cutting plane bound and picks a different next multiplier~$\lambda^3$.}

    \end{subfigure}
    \caption{After evaluating two multipliers $\lambda^1$. $\lambda^2$, we obtain solutions $X^1$ and $X^2$.
        Their implied inequalities reside on the convex shape, which is implicitly given and not known. The two generated cuts provide a 
        lower bound for $\mu^\ast$ given by $\mu^2$ (left). 
        An additional \emph{heuristic cut} can be created by arbitrary solutions (right).
        If it cuts off the next multiplier candidate $\mu^2$, we improve the gap without evaluating $\LR(\lambda)$. This cut will most likely not be a facet of the unknown boundary.}
    \label{fig:intuition}
    \end{figure}

In this section, we enhance the previously outlined general methodology by introducing a new strategy for identifying violated inequalities for the relaxation. As we will demonstrate, this enhancement not only improves the algorithm's performance but also addresses potential convergence issues. This improvement can have a substantial impact on practical outcomes, leading to a significantly better objective, which translates into improved long-term profitability.

 The original Cutting Plane procedure obtains cuts by solving~$\LR(\lambda)$. In our setting this involves heavy usage of (distributed) computational resources because determining objectives $f_i(x_i)$ still is costly -- each models the long-term profit for one article for a season and includes the determination of prices and sales in various countries for all weeks of the season. Thus, we investigate heuristic methods to obtain valid cuts with less overhead, i.e., given~$\lambda^j$ and~$\mu^j$ we would like to find $x \in \mathcal{X}$ such that Inequality~\eqref{ineq:cut} is violated for the current iterations multipliers. Figure~\ref{fig:intuition} illustrates this. The generation of such heuristic cuts has to trade-off speed vs. potential much stronger improvements by computing $\LR(\lambda^j)$. Furthermore, changing the strategy has a down-stream impact as the evaluated solutions $X^k$ are used to construct the final solution, see the last part of the previous section.

The core of our idea is to efficiently combine past found solutions~$X^1, \dots, X^j$ for which $f_i(X^k_i)$ is already computed for every iteration $k$ and article $i$. 
In that case we can find a new solution $X^{j+1}$ which is a combination of past solutions such that for all $i$ there is a $k \in {0,\dots,j}$ with $ X_i^{j+1} = X^k_i$.
Importantly, $X^{j+1}$ is a valid solution in $\mathcal{X}$, but not necessarily optimal for $\LR(\lambda, x)$ under any Lagrangian multiplier $\lambda$. Thus, there is no guarantee for finding a cut that maximizes the violation or even finding a violated cut at all.

Several strategies can be thought of to obtain $X^{j+1}$. We investigate the following three:
\begin{itemize}
\item Random heuristic: For each article $i$ sample $X_i^k$ uniformly from $k \le j$.
This gives us a baseline for comparing the other two strategies.
\item Maximum violation heuristic: Restrict the separation problem solved, i.e., the computation of~$\LR(\lambda^j)$ to already evaluated price strategies
\begin{equation} \label{eq:max_vio_heuristic}
 \sum_{i \in N} \max_{x_i\in \{X_i^k : k \le j\}} f_i(x_i)  - \lambda^\top A_{i} x_i.
\end{equation}
Since past evaluations of $f_i(X^k_i)$ as well as the contribution of each article to the linking constraints $A_{i} X^k_i$ are stored, this can be done efficiently by reweighing the violations and sorting.
\item Feasibility heuristic: Execute the primal heuristic from the previous section, using the current set of solutions (obtained from the preceding iterations). This process strives to derive a cut from a feasible and almost optimal solution. This approach is the most computationally involved among the three strategies discussed. Moreover, unlike the other two heuristics, repeatedly applying this strategy in sequence is not advantageous because it does not consider changes of the Lagrangian multipliers.
\end{itemize}

Note, that it does not suffice to only rely on these heuristic cuts. When the heuristic cutting process stalls, it can be beneficial to solve $\LR(\lambda)$, since this might help to get out of a "local optimum" when the cut information obtainable from $X^1, \dots, X^j$ is already or almost fully utilized. Finally, if no violated heuristic cut can be found, this might also be because the current solution value $\mu_j$ is in fact optimal for Problem~\eqref{prob:cp}, which, however, can only be proven by solving $\LR(\lambda)$.  

To make a trade-off between adding heuristic cuts and solving the Lagrangian relaxation we set a limit on the number of heuristic cut rounds before solving the Lagrangian relaxation again. Further, we employ two criteria to evaluate the usefulness of the latest heuristic cut and switched back to the exact separation of Inequality~\eqref{ineq:cut} if:
\begin{enumerate}
    \item The cut did not change the Lagrangian multipliers that minimize the cutting plane problem.
    \item The efficacy, a common measure for the quality of cuts (see, e.g.,~\cite{turner2023cutting}), falls below a defined threshold $\mbox{tol}_e$. In iteration $j$ the efficacy of the cut given by $X^{j}$ is given as
    \begin{equation}
     e(\lambda^j,\mu^j, X^j) = \frac{LR(\lambda^j, X^j) - \mu^j}{\| \lambda^j \|}.
    \end{equation}
    Figuratively, we stop if the generated cut is not separating $\mu^j$ for multiplier $\lambda^j$ deep enough.
\end{enumerate}

\begin{algorithm}[htb]
\caption{Extended cutting plane algorithm with heuristic cuts loop}
\label{alg:algorithm}
%\textbf{Input}:\\
\textbf{Parameter}: $\bar{n}$, $\bar{m}$, $\mbox{tol}_{\mu}$, $\mbox{tol}_{\nu}$, $\mbox{tol}_{e}$
%\textbf{Output}: 
\begin{algorithmic}[1] %[1] enables line numbers
\STATE Let $j\gets 0$, $\lambda^j \gets 0$.

\FOR{$n \gets 1$ to $\bar{n}$ }  \label{line:outer_loop}
    \STATE Compute $X^j$ as optimal solution of $LR(\lambda^j)$ and add the corresponding cut to Problem~\eqref{prob:cp_approx}
    \STATE Compute $\mu^j$, $\lambda^j$ as optimal solutions to Problem~\eqref{prob:cp_approx}
    %\STATE $\nu^n \gets \min_{k \le j} \LR(\lambda^k)$
    \STATE $j \gets j + 1$
    \FOR{$m \gets 1$ to $\bar{m}$} \label{line:inner_loop}        
        \STATE Compute a heuristic solution $X^j \in \mathcal{X}$ and add the corresponding cut to Problem~\eqref{prob:cp_approx}
        \STATE Compute $\mu^j$, $\lambda^j$ as optimal solutions to Problem~\eqref{prob:cp_approx}
        \IF {$\frac{\min_{k \le j} \LR(\lambda^k) - \mu^j}{\mu^j} < \mbox{tol}_{\mu}$ \OR  $\lambda^j = \lambda^{j-1}$ 
        % \OR  $\LR(\lambda^j, X^j) < (1-\mbox{tol}_{\nu}) \nu^n $ 
        \OR $e(\lambda^j,\mu^j, X^j) < \mbox{tol}_e$} \label{line:stop_inner_loop} 
        \STATE break
        \ENDIF
        \STATE $j \gets j + 1$
\ENDFOR
        \IF {$\frac{\min_{k \le j} \LR(\lambda^k) - \mu^j}{\mu^j} < \mbox{tol}_{\mu}$}   \label{line:stop_outer_loop}
        \STATE break
        \ENDIF
\ENDFOR
\end{algorithmic}
\label{alg:cuttingplane}
\end{algorithm}

Using the described heuristic with these stopping criteria we arrive at a modified algorithm for the cutting plane based Lagrangian descent, which is outlined in Algorithm~\ref{alg:cuttingplane} for the maximum violation heuristic.

In its outer loop starting in Line~\ref{line:outer_loop}, the algorithm performs up to $\bar{n}$ exact evaluations of multiplier $\lambda$.
Once the gap between dual bound and relaxed primal bound is within a configured tolerance, Line~\ref{line:stop_outer_loop} stops the algorithm.

After each multiplier evaluation, Line~\ref{line:inner_loop} adds up to $\bar{m}$ heuristic cuts to speed up the convergence. 
If an additional cut generated does not meet the efficacy tolerance $\mbox{tol}_e$, 
the resulting multipliers are not changing or the overall gap $\mbox{tol}_{\mu}$ is reached, Line~\ref{line:stop_inner_loop} ends the heuristic cut generation.

We did not perform intense parameter tuning, the chosen values are mostly ad-hoc with some sensitivity tests.
We limit the algorithm to $\bar{n}=10$ iterations. Up to $\bar{m}=100$ additional cuts could be generated. As tolerances, we ran our experiences with $\mbox{tol}_e = 1$ for the efficacy and $\mbox{tol}_{\mu} = 1\mathrm{e}{-6}$ for the dual bound gap. 

\subsection{Comparison with a disaggregated Cutting Plane formulation}\label{sec:disagg}
The cut heuristics are based upon the special characteristic of the markdown pricing problem: It is decomposable, i.e., the feasible set is the Cartesian product of the individual article level feasible sets, $\mathcal{X} = \mathcal{X}_1 \times \dots \times \mathcal{X}_n$. This characteristic actually allows to also obtain a so-called \emph{disaggregated}~\cite{guignard2003, frangioni2005} or \emph{multi-cut}~\cite{mitra2016} formulation 
\begin{equation}\label{prob:cp_approx_disaggregated}
    \begin{aligned}
     \mu^j, \lambda^j, \nu^j  =  \argmin_{\lambda \ge 0, \mu, \nu} \;&\mu  + \lambda^\top b 
        \\~\text{s.t.}~ & \mu  \ge \sum_{i \in N} \nu_i,\\
       &\nu_i\ge f_i(X^k_i) - \lambda^\top A X_i^k \quad\forall i \in N ~\forall k \le j.
        \end{aligned}
\end{equation}
Compare this with the formulation to obtain a bound after applying the maximum violation cut heuristic~$t$ times:
\begin{equation}\label{prob:cp_approx_cut_heuristic}
    \begin{aligned}
    \mu^j, \lambda^j  = \argmin_{\lambda \ge 0, \mu} \mu  + \lambda^\top b 
        ~\text{s.t.}~ \mu& \ge \sum_{i \in N} f_i(X^k_i) - \lambda^\top A X^k,\\
        \mu & \ge \sum_{i \in N} f_i(X^{k^s}_i)  - \lambda^\top  A X_i^{k^s} ~\forall 1 \le s < t,
        \end{aligned}
\end{equation}
where $k^s$ indicates the maximizing indices in Formula~\ref{eq:max_vio_heuristic} in round $s$ of applying the maximum violation heuristic.

One can show that the following strength guarantees on~$\mu^j$ from the different optimization problems holds: The disaggregated formulation~\eqref{prob:cp_approx_disaggregated} is at least as strong as the aggregated formulation with heuristic cuts~\eqref{prob:cp_approx_cut_heuristic}, which in turn is at least as strong as the aggregated formulation~\eqref{prob:cp_approx}. The first follows as the heuristic cutting planes are implied from the extended formulation of~\eqref{prob:cp_approx_disaggregated} via the summation of individual constraints. The latter follows, since~\eqref{prob:cp_approx} contains fewer constraints.

However, the aggregated formulation~\eqref{prob:cp_approx} contains only~$j+1$ variables and~$j$ constraints compared to~$j+s$ constraints with heuristic cuts and approximately $N+j$ variables and $N j$ constraints of the disaggregated formulation~\eqref{prob:cp_approx_disaggregated}. It turns out in our experiments, see Section~\eqref{sec:exp}, that solving~\eqref{prob:cp_approx_disaggregated} is not possible in practice due or very time consuming to its size. Thus using disaggregated cutting plane is not possible in our case. A compromise, which randomly partitions articles into~$M$ groups $N = N_1 \cup \dots \cup N_M $ to form the \emph{partially aggregated} formulation
    \begin{equation}\label{prob:cp_approx_partially_aggregated}
    \begin{aligned}
    \mu^j, \lambda^j, \nu^j  =  \argmin_{\lambda \ge 0, \mu, \nu} \;&\mu  + \lambda^\top b \\
        ~\text{s.t.}\; & \mu  \ge \sum_{1 \le \ell \le M} \nu_\ell,\\
       & \nu_\ell\ge \sum_{i \in N_\ell} f_i(X^k_i) - \lambda^\top A X_i^k  \quad\forall 1 \le \ell \le M ~ \forall k \le j
        \end{aligned}
    \end{equation}
is therefore evaluated below.

One can furthermore show that after a finite (but worst-case $j^N$) number of cut heuristics application, the computed bound is equal to the bound from the disaggregated formulation. To see this, note that after projecting~$\nu$ away in the disaggregated formulation, one obtains constraints of the form 
\begin{equation}\label{prob:heuristic_apporches_disagg}
    \mu  \ge \sum_{i \in N} f_i(X^{k_i}_i) - \lambda^\top A X_i^{k_i} \quad\forall i \in N ~\forall k_1, \dots, k_N, k_i \le j.
\end{equation}
Exactly this family of constraints is separated in the maximum violation cut heuristic. 

\section{Experiments and Real-world Impact}\label{sec:exp}
In the following section we illustrate the practical relevance of the methods described previously for the use case of the large-scale markdown price optimization system at Zalando. We report results here on a subset of the assortment for around 500,000 articles (depending on the exact time of the season), and for 14 countries, each with different discounts on weekly granularity until the season end of each article. Updated pricing recommendations are delivered at least twice a week, with 48 linking sales-weighted discount rate (sDR) constraints submitted as targets (24 lower and upper bounds each). To integrate with Zalando's business processes, the turn-around time between the setting of targets and the delivery of results should usually not exceed 3 hours.

Although considerable work has been put into improving the speed at which we can solve the article problems for one set of Lagrangian multipliers since we last reported on it~\citep{Li2022}, the time to solve the Lagrangian relaxation once remains at 150-700 CPU hours (using AWS C5 instances, exact time depending on the expected lifetime of the article). Even after parallelization over very large clusters the time to complete one iteration is therefore 3-12 minutes. Given additional overhead and because of both time and budgetary constraints we are limited to 10-15 iterations, a number that is far below the number of iterations one might usually expect for a problem with this number of constraints~\cite{guignard2003}.

In the experiments, we evaluate first which of the three cut generating strategies appears to be most effective, then the influence of the quality of the cuts on the dual bound. We go on to demonstrate that the cut heuristics are faster and easier to tune compared to partial aggregation strategies. Finally, we show how the associated commercial key performance indicators (KPI) benefit from the approach presented here. 

\subsubsection*{Comparison of cut generation heuristic strategies}
\begin{figure}[htb]
\includegraphics[width=13cm, keepaspectratio]{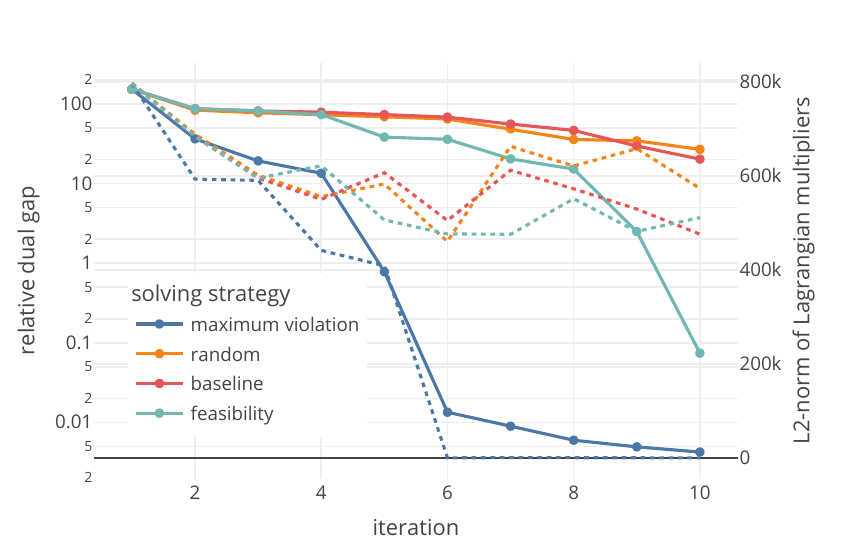}
\caption{Solid lines and left y axis: relative dual gap of the cutting plane problem over number of outer loop cutting plane iterations. Dashed lines and right axis: L2 norm of Lagrangian multipliers $\lambda^{j+1}$ selected in the cutting plane problem. Each iteration is one new solution of the Lagrangian dual. Baseline uses no heuristic to add cuts to the cutting plane problem. For details on the heuristics see Section~\ref{sec:heuristics}.}
\label{figure:allmethods}
\end{figure}

To explore whether the heuristics described above help improve the quality of solutions, we first tested the relative effectiveness of the different strategies. Fig.~\ref{figure:allmethods} illustrates the performance via a case study of a problem with problematic convergence in the baseline (no-heuristic) strategy; other instances behave similarly. We find that the maximum violation heuristic is most effective, closing the relative dual gap, $d_j = \frac{\min_{k \le j} \LR(\lambda^k) - \mu^j}{\min_{k \le j} \LR(\lambda^k)}$,
to less than 2\% after six outer loop iterations and 0.4\% in the 10th iteration for a particularly hard problem  which did not converge after 10 iterations in the baseline (with a relative gap much larger than 100\%).  We observed, however, that during the first few iterations, the relaxed primal bound hardly improved, and our imposed efficacy criteria would interrupt the heuristic cut generation early and continue with a full Lagrangian step.

Applying the feasibility heuristic is also effective in generating useful additional cuts, but less so. For that heuristic, the gap only drops below 100\% in the 10th iteration, with a final gap of 7.4\%. Random cuts are not effective, with the results being virtually indistinguishable from not applying any heuristic. We have observed similar results on other problem instances and concluded that the maximum violation heuristic strategy is performing the best. Thus we selected it as the most promising candidate for production and further results focus only on the maximum violation heuristic.

\subsubsection*{Comparing the maximum violation heuristic with partial aggregation formulations}

As discussed earlier, solving the Lagrangian relaxation quickly is of utmost importance to us. Therefore, we wanted to compare directly the effectiveness of the maximum violation cut heuristic and compare it to the previously reported partial aggregation formulation, see Section~\ref{sec:disagg} and ~\cite{guignard2003, frangioni2005,mitra2016}. Importantly, we want to ensure that the maximum violation strategy produces as good or better bounds as the partial aggregation formulation would in similar time.

In order to compare these strategies, we first look at their performance in isolation. We take dual problems that were previously identified as "difficult" as they required a large number of dual solutions (iterations) for the dual gap to close. For each of those problems we selected a set of dual solutions and applied one of two strategies:
\begin{itemize}
\item For the maximum violation heuristic, we take the dual solutions and repeatedly apply the heuristic. Each repetition reduces the dual gap, at the cost of time. We can therefore determine the dual gap as a function of time.
\item For the partially aggregated problem, we solve the optimization problem for different group sizes, starting from a single group (creating the equivalent of baseline cutting plane problem) up to a fully disaggregated problem. The assignment of articles into groups is done randomly. 

Larger group sizes will lead to larger problems and thus slower solving times, but also usually result in tighter gaps. Hence, we can visualize the gap as a function of elapsed time.
\end{itemize}

\begin{figure}[htb]
\centering
\includegraphics[width=10cm, keepaspectratio]{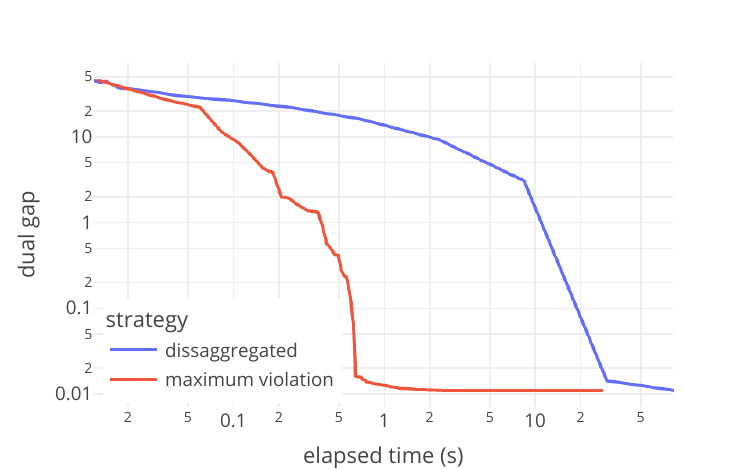}
\caption{Example of applying the partial aggregation model or the maximum violation strategy to a fixed set of dual solutions. For maximum violation (red), the heuristic is applied repeatedly, and the observed dual gap $d_j$ and elapsed time are noted. For the disaggregated formulation (blue), different levels of aggregation are solved independently, and the elapsed solving time and dual gap is noted.}
\label{figure:gap_line_plot}
\end{figure}

This is equivalent to taking a set of solutions $X^j$ as described in Algorithm \ref{alg:algorithm} and only applying steps 6 through 13. An example of a single problem is shown in Fig.~\ref{figure:gap_line_plot}. We can see that the maximum violation strategy closes the dual gap much quicker than the partially aggregated strategy for this example. 

To test if this generalizes, we collected a set of 86 other "difficult" dual problems. For all these instances, we were able to solve the fully disaggregated problem and thus obtained the best primal bound achievable given the parameters. As discussed in Section~\ref{sec:disagg}, the primal bounds provided by the disaggregated problem and by the maximum violation heuristic will eventually reach the same value. This is exactly what we see in our experiments: The bound provided by running the maximum violation heuristic until convergence was in all cases at the level of numeric precision of the best primal bound.

To compare the two strategies, we compute the time needed to reach a given gap to the best primal bound. For the heuristic, it is the time needed to achieve this gap by repeatedly applying the heuristic. For the partially aggregated strategy, we individually solve again the problem for different levels of aggregation and take the minimum
of the solving times for the aggregations leading to a primal bound within the gap.

Looking at the geometric means across the tested samples (Tab.~\ref{table:speedup}) we can see that the maximum violations strategy clearly outperforms the partially aggregated strategy in speed when solving to small gap sizes (0.1\% and 0.01\%), with the advantage increasing as the target gap gets smaller. For a gap of 0.01\%, the geometric mean of solving times for the partially aggregated strategy is 48 s, while the maximum violation strategy takes less than 2 s. Interestingly, for a larger target gap of~1\%, the geometric mean of solving times is nearly identical (0.12 s vs 0.11 s), but the geometric standard deviation is much larger for the partially aggregated problem.

\begin{table}[htb]
\centering
\begin{tabular}{lrr}
\toprule
 gap to best bound & partially aggregated & maximum violation \\
\midrule
0.1\% & 0.12 $\pm$ 18.43 & 0.11 $\pm$ 4.12 \\
0.01\% & 4.27 $\pm$ 10.24 & 0.59 $\pm$ 3.45 \\
0.001\% & 47.98 $\pm$ 3.35 & 1.69 $\pm$ 2.90 \\
\bottomrule
\end{tabular}
\caption{Geometric mean $\pm$ geometric standard deviation of solving time (s) for different gaps to best primal bound across 86 examples.}
\label{table:speedup}
\end{table}

\begin{figure}[htb]
\centering
\includegraphics[width=13cm, keepaspectratio]{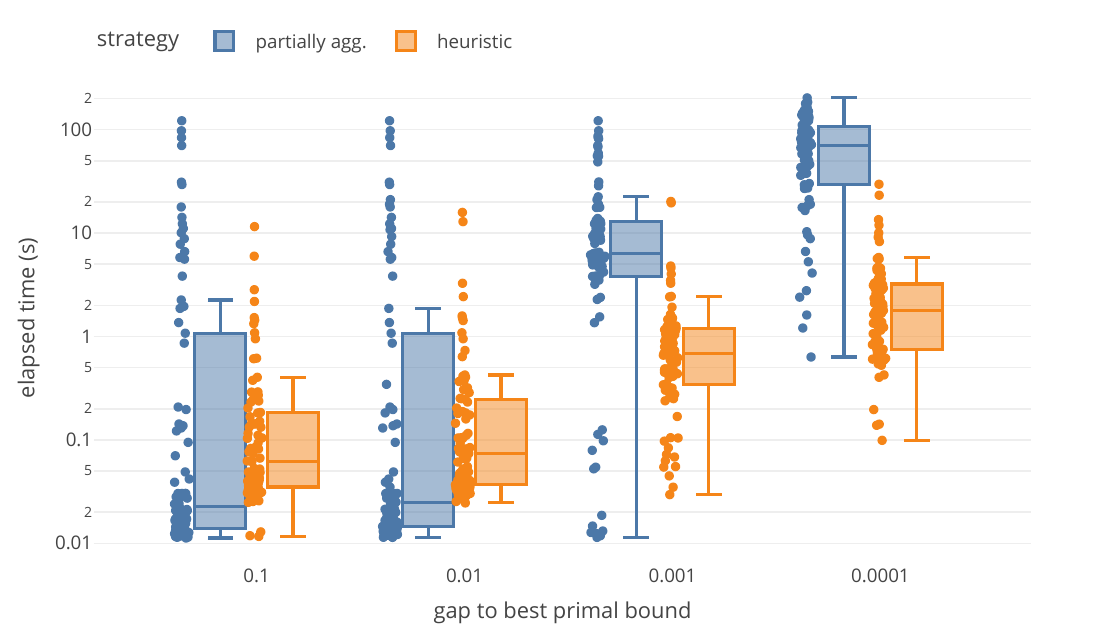}
\caption{Elapsed time to reach different gaps to best dual bound. For 86 different initial problems (a set of solutions to the relaxed dual problems, iterations from the Lagrangian decomposition) the partial aggregation and heuristic strategies were applied to find primal bounds. For partial aggregation, different group sizes were tested and used to determine the minimum time to gap. For the maximum violation heuristic, the heuristic was applied repeatedly, and the cumulative time to reach the target gap was recorded.}
\label{figure:direct_gap_size}
\end{figure}

The reason for this difference can be seen in the distribution of solving times for different gaps (Fig.~\ref{figure:direct_gap_size}). For larger target gaps (10\% or 1\%), the partially aggregated problem often solves them faster to that gap, with a lower median solve time, in cases where a partially aggregated problem with few groups was sufficient to solve to that target gap. At the same time, the longest solving times to the 10\% and 1\% gaps are much more extreme in the partially aggregated problems, sometimes exceeding 100 s, even though all problems took less than 20 s for the maximum violation heuristic. This split in the partially aggregated problems, where some problems are solved fast, while others are very slow, is also seen in the 0.1\% gap target (with more problems taking longer than for the maximum violation heuristic). For the 0.01\% target gap  we see the heuristic very consistently outperforming.

All-in-all, this leaves us with several clear advantages for the maximum violation heuristic over the partial aggregation strategy in our setup:
\begin{itemize}
\item The maximum violation heuristic reaches tight primal bounds faster than the partial aggregation strategy. Particularly for tight bounds, the heuristic is more than 20 times faster (in geometric mean).
\item The maximum violation heuristic produces more consistent bounds at a given time. This consistency makes it easier to tune, as one can easier select how much time will usually result in "tight enough" bounds.
\item Because the primal bound can be checked on every application of the maximum violation heuristic, it is easier to define stopping criteria.
\end{itemize}

\begin{figure}[htb]
\includegraphics[width=13cm, keepaspectratio]{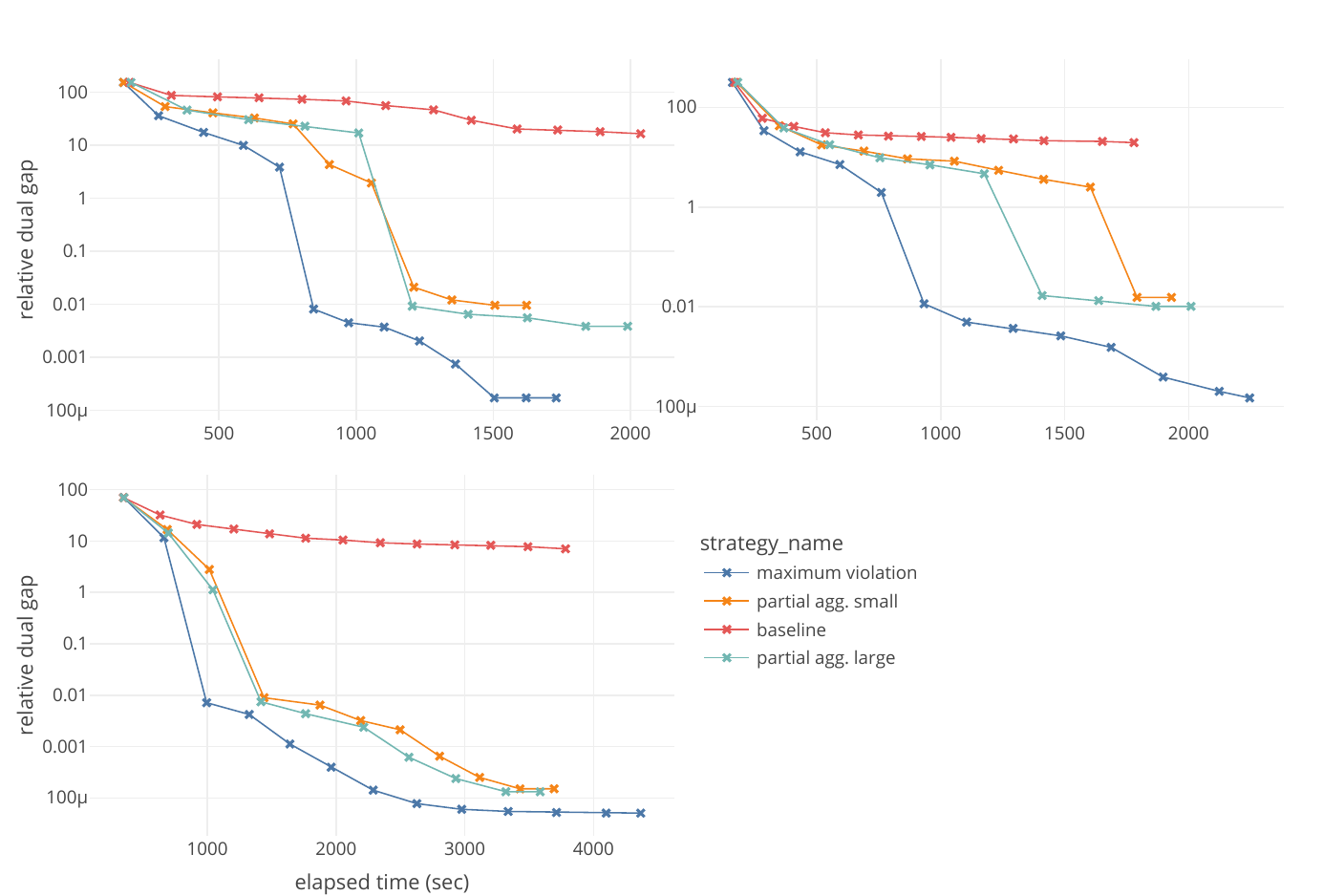}
\caption{Full Lagrangian descent using different solving strategies. For 3 different problems (different weeks), the maximum violation strategy was applied (blue, 100 times per solving of the LR), and compared to baseline (red, baseline cutting plane), partially aggregated small (orange, creating 2000 random groups), and partially aggregated large (cyan, using 5000 groups). For each problem the dual gap $d_j$ was noted after each solving of the Lagrangian relaxation.}
\label{figure:gap_size_diff}
\end{figure}

Still, we want to compare how well these strategies perform when applying them to the complete solving of dual problems (that is, the full Algorithm~\ref{alg:algorithm}). Here, we need to balance between the time spend obtaining dual solutions (solutions to $LR(\lambda^j)$) and getting tight primal bounds and good corresponding Lagrangian multipliers $\lambda^j$. As we have seen in Fig.~\ref{figure:direct_gap_size}, this is particularly challenging for the partial aggregation strategy, as the quality of the bound for a given aggregation level can vary drastically between problems. Therefore, based on those results, we decided to try two different partial aggregation levels: the sub-problems are aggregated to either groups of 2000 (small) or 5000 (large) articles. Along with the maximum violation heuristic (adding 100 heuristic cuts on each multiplier evaluation) and the baseline (baseline cutting plane without either partial aggregation or cut heuristic), these strategies are applied to 3 problems that were previously found to be difficult to solve. The dual gaps for these problems are shown against the elapsed time in Fig.~\ref{figure:gap_size_diff}. 

We see that the maximum violation heuristic (blue) provides a marked improvement over baseline (red), which fails to close the dual gap in all 3 examples. Interestingly, the maximum violation heuristic also outperforms the partially aggregated formulation for both aggregation levels. While both partially aggregated problems do provide reasonable gaps on the problems (with gaps at or below 0.01\%), the maximum violation strategy reaches smaller gaps (around 0.0001\%), and closes the gap faster. The maximum violation heuristic also takes less time per multiplier evaluation (outer loop in Algorithm~\ref{alg:algorithm}) compared to the partially aggregated problems. Given that reaching very tight bounds consistently would require more than 100 s per multiplier evaluation, it is therefore unlikely that any further increasing of the group size would change the advantage we see for the maximum violation heuristic.

Based on these findings, i.e., quicker convergence, higher consistency and better tunability, we decided to select the maximum violation heuristic over the partial aggregation strategy for further experimentation and the real application.

\subsubsection*{Impact on commercial KPIs}

For a fair comparison of the overall improvement that applying the maximum violation heuristic gives in a time-constrained setting, we went on to compare final solutions (after application of the primal heuristic described in Section~\ref{sec:heuristics}) both with the maximum violation heuristic cut strategy and without a heuristic cut strategy. 

Since the use of the maximum violation heuristic cut strategy does not significantly increase run times, both strategies would take approximately equal time to solve and the results should be directly comparable. In this final experiment we want to investigate whether an improvement in solving the dual problem also leads to an improvement of the primal solution generated after applying our primal heuristic by solving problem~\eqref{prob:primal_heuristic}.

Testing the improvement for 8 different dates we can see that through the application of the maximum violation heuristic in forward running simulations we achieve an uplift in all critical commercial KPIs: long-term profit (LTP, the primal objective in our formulation), first week gross merchandise volume (GMV) and the first week profit contributor 2 (PC2, which is a company accounting metric, that measures revenue minus the combined cost of sales and fulfilment). For the selected dates, we find an average weekly improvement of more than 3M€ in LTP and GMV and a 0.9M€ improvement in PC2 (Table \ref{table:LTPimprovement}). Our simulation indicates that using our suggested strategy could lead to a long term profit-increase of about 27M€ in eight weeks, which projects to 175M€ per year.

\begin{table}[htb]
\centering
\begin{tabular}{lrrr}
\toprule
 & LTP & GMV & PC2 \\
Experiment &  &  &  \\
\midrule
0 & 0.18M€ & 0.31M€ & 0.09M€ \\
1 & 0.55M€ & 0.87M€ & 0.23M€ \\
2 & 3.30M€ & 6.02M€ & 1.42M€ \\
3 & 3.63M€ & 6.07M€ & 1.54M€ \\
4 & 3.65M€ & 1.38M€ & 0.75M€ \\
5 & 3.97M€ & 4.74M€ & 1.28M€ \\
6 & 5.71M€ & 0.00M€ & 0.00M€ \\
7 & 6.20M€ & 7.25M€ & 2.09M€ \\
\midrule
avg &  3.40M€ &  3.33M€ & 0.93M€ \\
sum & 27.19M€ &  26.64M€ & 7.40M€ \\
\bottomrule
\end{tabular}
\caption{Simulated improvement (in M€) of the primal solution in key KPIs for 8 start dates throughout the season. Improvements are the difference between the primal solution generated without the maximum violation heuristic, and the solution generated by applying the maximum violation heuristic 100 times per iteration. The last two rows show the average and the sum of the eight weeks}
\label{table:LTPimprovement}
\end{table}

The significant enhancements observed in our simulations, combined with the minimal risk of solution deterioration (given that the heuristic cuts added are still valid, just potentially less effective), have convinced us and the responsible stakeholders to implement the maximum violation heuristic in Zalando's production environment. 

Due to parallel, unrelated AB tests, we were unable test the improvements through direct comparisons, and instead opted to follow the improvements using causal impact analysis~\citep{causalimpact}. In causal impact analysis, we select predictors that can be used to train a Bayesian structural time-series model to forecast the variables of interest. Assuming that the predictors are unaffected by the treatment we can use them to forecast what would have happened if no intervention had taken place. See also \cite{mehrotra2020price} for another usage of causal impact analysis to estimate treatment effects without relying on AB testing.

Since long term profit (the immediate objective in our markdown optimization) is not directly measurable over short observation horizons, we instead opted to focus on the improvements in PC2 and GMV we observed in our forward running simulations. We chose to use 25 weeks prior to treatment (the application of the maximum violation heuristic) for training of the time-series model and monitored the effect of treatment for 10 weeks. As covariates for the time-series model we chose the respective predictions of the markdown optimization model without the maximum violation heuristic for PC2 and GMV. Since these covariates are obtained by computing the old optimization without heuristic cuts and are therefore unaffected by the treatment (the use of the maximum violation heuristic in the model), this gives us a good estimator of the observed effect.

\begin{figure}[htbp]
\includegraphics[width=13cm, keepaspectratio]{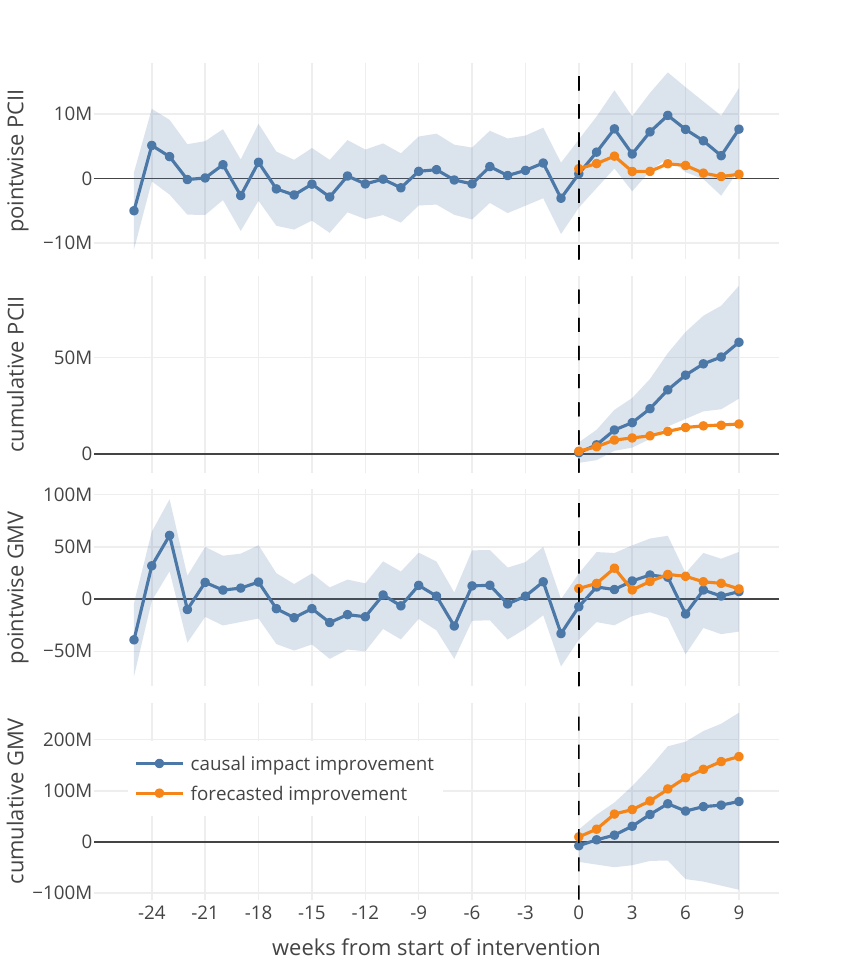}
\caption{Causal impact analysis of the switch to the maximum violation heuristic. Analysis was performed using 25 weeks prior and 10 weeks post switching the heuristic online, using an optimizer without the heuristic as predictor and an observed PC2 and GMV as variable. We show the uplift in percent of the average weekly performance for GMV and PC2. Blue are the causal impact effects (observed - Bayesian structural time-series forecast), orange are the effects comparing optimization results without heuristic to those with heuristic.}
\label{figure:causalimpact}
\end{figure}

In the causal impact model (Fig.~\ref{figure:causalimpact}) we see an improvement in PC2 in the treatment period (p=0.0001, observed weekly effect: 5.8M€, expected effect 1.56M€) and an effect on GMV that is not significant, but consistent with expectations (p=0.18, observed weekly effect: 7.9M€, expected effect: 16.7M€) These findings indicate that the improvement computed in the forward running simulations are actually observable in practice, and that applying the maximum violation heuristic leads to tangible real world improvements.

\section{Conclusion}\label{sec:concl}

In this paper, we addressed the practical challenges encountered in implementing large-scale markdown strategies using a predict-then-optimize framework. These challenges primarily revolve around the stability and efficiency of identifying close-to-optimal solutions for huge mixed integer programs. We introduced novel heuristics that enhance the performance of existing methods, facilitating the rapid identification of high-quality solutions through a Lagrangian relaxation approach. An important component was the implementation of a heuristic cut generation scheme, based on a maximum-violation measure, and its seamless integration within an existing exact separation procedure. We demonstrated the practical impact of our methodology with empirical evidence from Zalando SE's pricing systems. Our approach not only accelerates the process of finding high-quality solutions but also leads to multi-million revenue and profit increases, underscoring its commercial relevance.

\bibliographystyle{elsarticle-harv} 

\bibliography{main.bib}
\clearpage

\appendix

\end{document}